\documentclass{ifacmtg}
\usepackage[dcucite]{harvard}

\usepackage[dvips]{graphicx}
\usepackage{makeidx}
\usepackage{amsmath}
\usepackage{amsfonts}

\newcommand {\alphavec}{\boldsymbol{\alpha}}

\newcommand {\xivec}{{\boldsymbol{\xi}}}

\newcommand {\nuvec}{\boldsymbol{\nu}}
\newcommand {\varevec}{\boldsymbol{\varepsilon}}

\newcommand {\thetavec}{{\boldsymbol{\theta}}}

\newfont{\pseudocode}{cmtt10}

\newtheorem{assume}{Assumption}

\newcommand{\bfx}{\mathbf{x}}
\newcommand{\bff}{\mathbf{f}}

\newcommand{\bfz}{\mathbf{z}}

\newcommand{\bfu}{\mathbf{u}}

\newcommand{\dpsi}{{\dot\psi}}
\newcommand{\pd}{{\partial}}

\newcommand{\Real}{\mathbb{R}}

\makeindex

\begin{document}



\runauthor{Ivan Tyukin, Denis Efimov, Cees van Leeuwen}
\begin{frontmatter}
\title{Adaptive Regulation to Invariant Sets}
\author[Wako-shi]{Ivan Tyukin}
\author[IPME]{Denis Efimov}
\author[Wako-shi]{Cees van Leeuwen}
\address[Wako-shi]{Laboratory for Perceptual Dynamics,
RIKEN Brain Science Institute, 2-1, Hirosawa, Wako-shi, Saitama,
Japan}
\address[IPME]{Institute for Problems of Mechanical Engineering, Laboratory for Control of Complex
Systems, V. O., Bolshoy  61, Saint-Petersburg, Russia}

\begin{abstract}
A new framework for adaptive regulation to invariant sets is
proposed.  Reaching the target dynamics (invariant set) is to be
ensured by state feedback while adaptation to parametric
uncertainties is provided by additional adaptation algorithm. We
show that for a sufficiently large class of nonlinear systems it
is possible to adaptively steer the system trajectories to the
desired non-equilibrium state without requiring  knowledge or
existence of a specific strict Lyapunov function.
\end{abstract}

\end{frontmatter}

 \vspace{5mm}

\section{Introduction}

Whether adaptive or non-adaptive solutions are sought in  control
theory, the problem is usually stated in terms of stabilization
problem of an equilibrium or tracking of a given reference signal.

In recent years, motivated by problems in physics and natural
sciences, slightly different demands came to the surface. Instead
of forcing a system to an arbitrary equilibrium one should search
for the natural motions in the system which satisfy the control
goal the most, and then transform these to the desired state by
gentle and small control efforts \cite{Kolesn,Fradkov2003}. One of
the successful examples is the result reported in the seminal
paper \cite{Ott_1990} with long-standing theoretical impact and
exciting practical applications \cite{Tziperman_1997}.

The problem with this and similar methods, however, in the context
of adaptive control is that for the given feedback one must know
the Lyapunov function ensuring asymptotic stability of the target
dynamics. This gives rise to another severe limitation of the
conventional Lyapunov-based methodology  -- the problem with
asymptotic behavior of adaptive systems. Roughly speaking, the
problem is as follows: while the specific Lyapunov function fits
very well the non-adaptive controller design (i. e. ensures that
solutions converge asymptotically to the desired state), it may
not guarantee the desired asymptotic in the adaptive case. The
reason for such is that the Lyapunov function itself is not
strict. The breakthrough in this problem has been reported in
\cite{Panteley_2002,Ortega2003}. The problem has been resolved for
equilibria that can be made asymptotically stable by state
feedback. Yet, non-equilibrium and non-asymptotically stable
dynamics were not addressed.

The problems of non-equilibrium control are gaining substantial
attention in the recent years, especially in the framework of
output regulation. In \cite{Byrnes_2003} a number of sufficient
and necessary conditions assuring existence of the solution to
this problem are proposed. Although {\it the internal model
principle} in the output regulation problem \cite{Byrnes_1997}
proves strong bindings between adaptive and output regulation
problems, historical and methodological differences in these
branches of the control theory do not always allow explicit
application of the results from one field to another. This
provides additional motivation to our current study in the context
of adaptation.

The contribution of our present paper is as follows. First, we aim
to formulate the problem of adaptive regulation to the desired
non-equilibrium dynamics. This dynamics should in principle be
invariant under the system flow. It also should poses certain
properties like boundedness of the trajectories and/or partial
stability \cite{Vorotnikov}. No asymptotic Lyapunov-like stability
conditions are to be imposed a-priori in order to escape the
burden of detectability. Second, under these assumptions we shall
be able to derive adaptation algorithms which are capable of
steering the system trajectories to the desired invariant set. In
order to do so we employ the recently developed  {\it adaptive
algorithms in finite form} \cite{t_fin_formsA&T2003}. These
algorithms guarantee improved performance and are capable of
handling nonlinear parametrization of the uncertainty
\cite{t_fin_forms_arch}. The main idea of this approach is to
introduce the desired invariant set into the system dynamics
(virtual adaptation algorithms) and then realize these algorithms
by means of the embedding technique proposed in
\cite{ECC_2003,ALCOSP_2004,t_fin_forms_arch}.

The paper is organized as follows: in Section 2 we provide
necessary notations and formulate the problem. Section 3 contains
the main results of the paper given in Theorem
\ref{theorem:adaptive_invariant}. The proof of the theorem is
provided in the subsequent subsections. Each of the subsection
substitutes the separate step of the construction. Subsection
\ref{proof:part_1} addresses design of the virtual algorithms,
Subsection \ref{proof:part_2} provides auxiliary system which is
necessary for the embedding, Subsection \ref{proof:part_3}
contains the main arguments of the proof. Section 4 concludes the
paper.

Throughout the paper we will use the following notations: symbol
$\bfx(t,\bfx_0,t_0)$ stands for the flow which maps
$\bfx_0\in\Real^{n},t_0,t\in \Real_+$ into $\bfx(t)$.  Function
$\nu: R_+\rightarrow R$ is said to belong to $L_2$ iff
$L_2(\nu)=\int_0^{\infty}\nu^2(\tau)d\tau<\infty$. The value
$\sqrt{L_2(\nu)}$ stands for the $L_2$ norm of $\nu(t)$. Function
$\nu: R_{+}\rightarrow R$ belongs to $L_\infty$ iff
$L_{\infty}(\nu)=\sup_{t\geq0}\|\nu(t)\|<\infty$, where
$\|\cdot\|$ is the Euclidean norm. The value of $L_{\infty}(\nu)$
stands for the $L_{\infty}$ norm of $\nu(t)$.

\section{Problem formulation}

\begin{defn} A point $p \in \Real^n$ is called an
$\omega$-limit point $\omega(\bfx(t,\bfx_0,t_0))$ of $\bfx_0\in
\Real^n$ if there exists a sequence $\{t_i\}$,
$t_i\rightarrow\infty$, such that $\bfx(t,\bfx_0,t_0)\rightarrow
p$. The set of all limit points $\omega(\bfx(t,\bfx_0,t_0))$ is
the $\omega$-limit set of $\bfx_0$.
\end{defn}
In order to specify explicitly in our notations which particular
flow is referred to in the notion of the $\omega$-limit set we use
notations $\omega_\bff(\bfx_0)$ (and $\bfx_\bff(t,\bfx_0,t)$) to
denote the $\omega$-limit set (and flow) of $\bfx_0$ in the
following system $\dot{\bfx}=\bff(\bfx), \ \bfx_0\in X\subset
\Real^{n}$. Symbol $\Omega_{\bff}(\bfx)$ denotes the union of all
$\omega_\bff(\bfx_0), \ \bfx_0\in X$. Throughout the paper we will
refer to  set $\Omega_{\bff}(\bfx)$ as the $\Omega_\bff$-limit set
(or simply $\Omega$-limit set if the corresponding flow is defined
from the context) of the system.

\begin{defn} Set $S\subset \Real^n$ is invariant (forward-invariant) under the flow
$\bfx_\bff(t,\bfx_0,t_0)$ iff $\bfx_\bff(t,\bfx_0,t_0)\in S$ for
any $\bfx_0\in S$ for all $t>t_0$.
\end{defn}

In our current study we consider the following class of systems:
\begin{equation}\label{system_affine}
\begin{split}
\dot{\bfx} & =\bff(\bfx)+G_u (\phi(\bfx)\thetavec+\bfu),
\\
\dot{\thetavec}&=S(\thetavec), \ \thetavec(t_0)\in \Theta\subset
\Real^{d}
\end{split}
\end{equation}
where $\bff:\Real^{n}\rightarrow \Real^n$,
$\phi:\Real^{n}\rightarrow \Real^{m\times d}$, are $C^0$-smooth
vector-fields, $G_u \in \Real^{n\times m}$,  $\thetavec$ is the
vector of unknown time-varying parameters, $S:
\Real^{d}\rightarrow \Real^{d}$, $S\in C^1$ is known, vector of
initial conditions ${\thetavec}(t_0)\in \Theta$, however, is
assumed to be unknown. Without loss of generality we assume that
$\Omega_S(\Theta)\subseteq \Theta$, and that $\Theta$ is bounded.
Our goal is to steer the state to the {\it target domain}:
\[
\Omega^{\ast}(\bfx)\subset \Real^{n}
\]

Let us introduce the following set of assumptions related to the
choice of  domain $\Omega^{\ast}(\bfx)$.

\begin{assume}\label{assume:target_boundedness} Set $\Omega^{\ast}(\bfx)\subset \Real^{n}$ is the
bounded and closed  set in $\Real^{n}$.
\end{assume}

\begin{assume}\label{asume:internal_model} There exists positive-definite matrix $H=H^T\in \Real^{d\times d}$, such that function $S:\Real^d\rightarrow \Real^d$ in
(\ref{system_affine}) satisfies the following inequality:
\[
H\frac{\pd S(\thetavec)}{\pd \thetavec}+\frac{\pd
S(\thetavec)}{\pd \thetavec}^{T}H\leq 0 \ \forall \
\thetavec\in\Real^d
\]
\end{assume}

\begin{assume}\label{assume:target_feasibility} For the given $\Omega^{\ast}(\bfx)$ and system
(\ref{system_affine}) there exists control function $\bfu_0(\bfx)$
such that
\[
G_u\bfu_0(\bfx)+\bff(\bfx)=\bff_0(\bfx)
\]
and, furthermore, for any $\bfx_0\in \Real^n$ the following holds:
$\Omega^{\ast}(\bfx)\subset \Omega_{\bff_0}(\bfx)$, where the flow
$\bfx(t,\bfx_0,t)$ is defined by
\begin{equation}\label{system_unperturbed}
\dot{\bfx}=\bff_0(\bfx)
\end{equation}
\end{assume}

Let us finally introduce two alternative hypotheses. The first
hypothesis is formulated in Assumptions
\ref{assume:target_invariance}, \ref{asseme:manifold_stability},
and \ref{assume:constraint}.  The second is given by Assumption
\ref{assume:iISS}.

\begin{assume}\label{assume:target_invariance} There exist functions  $\psi(\bfx): \ \Real^{n}\rightarrow
\Real$, $\varphi:\Real^{n}\rightarrow \Real$, and induced by
function $\psi(\bfx)$ set:
\[
\Omega_{\psi}=\{\bfx\in\Real^{n}| \ \bfx: \varphi(\psi(\bfx))=0 \}
\]
such that the following holds
$\Omega^{\ast}\subseteq\Omega_{\bff_0}(\Omega_\psi)$, i. e.
$\Omega^{\ast}(\bfx)$ is the largest invariant set of
(\ref{system_unperturbed}) in $\Omega_\psi$.
\end{assume}

\begin{assume}\label{asseme:manifold_stability} For the given function $\psi(\bfx): \Real^{n}\rightarrow
\Real$, $\psi(\bfx)\in C^1$ and vector field  $\bff_0(\bfx)$
defined in (\ref{system_unperturbed}) there exists function
$\beta(\bfx):\Real^{n}\rightarrow \Real_+$ such that $\beta(\bfx)$
is separated from zero and satisfies the following equality:
\begin{equation}\label{varphi}
\begin{split}
& \psi\frac{\pd \psi(\bfx)}{\pd \bfx}\bff_0(\bfx)\leq
-\beta(\bfx)\varphi(\psi)\psi,  \\
& \int_{0}^\psi \varphi(\sigma)d\sigma \geq 0, \
\lim_{\psi\rightarrow\infty} \int_{0}^\psi
\varphi(\sigma)d\sigma=\infty
\end{split}
\end{equation}
\end{assume}

\begin{assume}\label{assume:constraint} For the given function
$\psi(\bfx):\Real^{n}\rightarrow \Real$, $\psi(\bfx)\in C^1$ the
following relation holds:
\[
\psi(\bfx(t))\in L_\infty \Rightarrow \bfx\in L_\infty
\]
\end{assume}

Notice that function $\psi(\bfx)$ in Assumptions
\ref{asseme:manifold_stability}, \ref{assume:constraint} should
not necessarily be the (positive) definite function. Function
$\varphi(\psi)\psi$ is also not required to be (positive)
definite.

\begin{assume}\label{assume:iISS} Consider system
(\ref{system_unperturbed}) with additive input
$\varevec_0(t):\Real\rightarrow \Real^n$, $\varevec_0(t)\in
C^{1}$:
\begin{equation}\label{system_unperturbed1}
\dot{\bfx}=\bff_0(\bfx)+\varevec_0(t), \ \varevec_0\in L_2
\end{equation}
System (\ref{system_unperturbed1}) has finite $L_2\rightarrow
L_\infty$ gain, and in addition $\Omega^{\ast}\subseteq
\Omega_{\bff_0}$.
\end{assume}

The main question of our current study is that wether or not it is
possible to design the adaptation algorithm $\hat{\thetavec}(t)$
for system (\ref{system_affine}) such that the feedback of the
following form
\[
\bfu(\bfx,\xivec)=\bfu(\bfx,\xivec,\hat{\thetavec}), \
\dot{\xivec}=\bff_\xi(\bfx,\xivec,\hat{\thetavec}), \ \xivec\in
\Real^k
\]
ensures boundedness of the trajectories in the closed loop system
and that $\bfx(t)\rightarrow\Omega^{\ast}$ as
$t\rightarrow\infty$.

\section{Main Results}

The main idea of our approach is two-fold. First,  we search for
the desired dynamics of the closed loop system with feedback
$\bfu(\bfx,\xivec,\hat{\thetavec})$ and yet unknown
$\hat{\thetavec}(t)$, $\xivec(t)$ which ensures desired properties
of the controlled system. These properties should allow us to show
that under specific conditions $\bfx(t)\rightarrow\Omega^{\ast}$
as $t\rightarrow\infty$. Derivative of function
$\hat{\thetavec}(t)$ with respect to $t$ at this stage can, in
principle, depend on unknown parameters $\thetavec$. Family of all
such desired subsystems is referred to as {\it virtual adaptation
algorithms}.

The second stage of our method is to render these algorithms into
computable and physically realizable form. In particular, these
realizations should neither rely on a-priory unknown parameters,
nor should they require measurements of the right-hand side of
(\ref{system_affine}) (i.e. derivatives).

In order to achieve this goal we invoke the {\it algorithms in
finite form} \cite{t_fin_formsA&T2003,t_fin_forms_arch}
(physically realizable and computable control) and the embedding
argument inctroduced in \cite{ECC_2003,ALCOSP_2004}. In general,
finite form realizations of virtual adaptation algorithms require
analytic solution of a partial differential equation known as {\it
explicit realization condition}. However, with the embedding
technique proposed in our earlier publications it is possible to
avoid this difficulty and derive adaptation schemes as a known and
well-defined function of $\bfx,t$. The main result of our current
study is formulated in Theorems \ref{theorem:adaptive_invariant}
and \ref{theorem:adaptive_invariant_integral} below.
\begin{thm}\label{theorem:adaptive_invariant} Let system (\ref{system_affine_aux}) be given and
Assumptions
\ref{assume:target_boundedness}--\ref{assume:constraint} hold.
Let, in addition, there exists $C^1$-smooth function
$\kappa(\bfx)$ such that the following estimate holds:
$\left\|\frac{\pd \psi(\bfx)}{\pd \bfx}\right\|\leq
|\kappa(\bfx)|$. Then there exists auxiliary system
\begin{equation}\label{systen_ext_gen}
\begin{split}
\dot{\xivec}&=\bff_{\xi}(\bfx,\xivec,\nuvec) \\
\dot{\nuvec}&=\bff_{\nu}(\bfx,\xivec,\nuvec), \ \xivec\in
\Real^{n}, \ \nuvec\in \Real^{d}
\end{split}
\end{equation}
control input
$\bfu(\bfx,\hat{\thetavec})=\bfu_0(\bfx)-\phi(\xivec)\hat{\thetavec}(t)$,
and adaptation algorithm
\begin{equation}\label{fifo_lin_invar_set}
\begin{split}
\hat{\thetavec}&=(H^{-1}\Psi(\xivec)\bfx+\hat{\thetavec}_I(t)),
\\
\Psi(\xivec)&=(\kappa^2(\xivec)+1)(G_u\phi(\xivec))^{T} \\
\dot{\hat{\thetavec}}_I&= S(\hat{\thetavec})-H^{-1}\frac{\pd
\Psi(\xivec)}{\pd \xivec} \bff_{\xi}(\bfx,\xivec,\nuvec)
\bfx- \\
&   H^{-1}\Psi(\xivec)\bff_0(\bfx)
\end{split}
\end{equation}
such that the following properties hold:

1) $\hat{\thetavec}(t),\bfx(t)\in L_\infty$

2) trajectories $\bfx(t)$ converge into the domain $\Omega^{\ast}$
as $t\rightarrow \infty$

3) if $G_u\phi(\xivec)$ is persistently exciting then
$\hat{\thetavec}(t,\hat{\thetavec}_0,t_0)$ asymptotically
converges to $\thetavec(t,\thetavec_0,t_0)$.
\end{thm}

\begin{thm}\label{theorem:adaptive_invariant_integral} Let system (\ref{system_affine_aux}) be given and
Assumptions
\ref{assume:target_boundedness}--\ref{assume:target_feasibility},
and \ref{assume:iISS} hold. Then there exist auxiliary system of
type (\ref{systen_ext_gen}), control input
$\bfu(\bfx,\hat{\thetavec})=\bfu_0(\bfx)-\phi(\xivec)\hat{\thetavec}(t)$
and adaptation algorithm (\ref{fifo_lin_invar_set}) with
$\kappa(\xivec)\equiv 0$ such that statements 1)--3) of Theorem
\ref{theorem:adaptive_invariant} hold.
\end{thm}

The proof of the theorems is given in the next subsections. In
subsection \ref{proof:part_1} we derive virtual adaptation
algorithms which satisfy in part the requirement of the theorem.
Subsection \ref{proof:part_2} introduces function $\xivec(t)$
satisfying the embedding assumption from
\cite{ECC_2003},\cite{t_fin_forms_arch}. In subsection
\ref{proof:part_3} we combine these results together and complete
the proofs.

\subsection{Design of Virtual Adaptive Algorithms}\label{proof:part_1}

Let us consider the following dynamic state feedback
$\bfu(\bfx,\hat{\thetavec})=\bfu_0(\bfx)-\phi(\xivec)\hat{\thetavec}(t)$.
This feedback renders system (\ref{system_affine}) into the
following form
\begin{equation}\label{system_affine_x}
\begin{split}
\dot{\bfx}&=\bff_0(\bfx)+G_u\phi(\xivec)(\thetavec-\hat{\thetavec}(t))+
\\ &  G_u(\phi(\bfx)-\phi(\xivec))\thetavec,
\end{split}
\end{equation}
Let us denote $G_u\phi(\bfx)=\alphavec(\bfx)$ and consider the
following auxiliary system
\begin{equation}\label{system_affine_aux}
\begin{split}
\dot{\bfx}&=\bff_0(\bfx)+\alphavec(\xivec)(\thetavec-\hat{\thetavec})+\varevec(t),
\\
\dot{\thetavec}&= S(\thetavec)\\
 \dot{\hat{\thetavec}}&=
S(\hat{\thetavec})+H^{-1}(\kappa^2(\xivec)+1)\alphavec(\xivec)^{T}\times \\
&  (\alphavec(\xivec)(\thetavec-\hat{\thetavec})+\varevec(t))\\
\kappa(\xivec)&:\Real^{n}\rightarrow \Real, \ \kappa\in C^1
\end{split}
\end{equation}

\begin{lem}[Virtual Adaptation Algorithm]\label{lemma:virtual} Let system (\ref{system_affine_aux}) be given and
Assumptions
\ref{assume:target_boundedness}--\ref{assume:target_feasibility},
\ref{assume:constraint} hold. Furthermore, let
$\kappa(\xivec(t))\varevec(t)\in L_2$, and $\varevec\in L_2$.

Then the following statements hold:

1) $\hat{\thetavec}(t)$ is bounded for every
$\theta(t_0)\in\Theta$, $\hat{\thetavec}(t_0)\in\Real^{d}$

2)
$\kappa(\xivec)\alphavec(\xivec)(\hat{\thetavec}(t)-\thetavec(t))$,
$\alphavec(\xivec)(\hat{\thetavec}(t)-\thetavec(t))\in L_2$

3) Let, in addition, $\left\|\frac{\pd \psi(\bfx)}{\pd
\bfx}\right\|\leq |\kappa(\bfx)|, \ \ \bfx-\xivec\in L_\infty$
then $\bfx\in L_\infty$

4) if, independently on the conditions of statement 3),
$\varevec(t)\equiv 0$ and the function $\alpha(\xivec)$ is
persistently exciting, i. e. there exist constants $\delta, T>0$
such that
$\int_{t}^{t+T}\alpha(\xivec(\tau))^{T}\alpha(\xivec(\tau))\geq
\delta I_{d}$ then trajectory $\hat{\thetavec}(t)$ converges to
the solution $\thetavec(t,\thetavec_0,t_0)$ exponentially fast.
\end{lem}
{\it Lemma \ref{lemma:virtual} proof.} Let us show that statements
1) and 2) hold. Consider the following positive-definite function:
\[
V_{\theta}(\thetavec,\hat{\thetavec},t)=\|\thetavec-\hat{\thetavec}\|^{2}_{H}+\epsilon=(\thetavec-\hat{\thetavec})^T
H
 (\thetavec-\hat{\thetavec}) +\epsilon,
\]
where
$\epsilon(t)=\frac{1}{2}\int_t^{\infty}(\kappa^2(\xivec(\tau))+1)\varevec^{T}(\tau)\varevec(\tau)d\tau\geq
0$. According to the lemma assumptions function
$\kappa(\xivec(t))\varevec(t)\in L_2$. This implies that
$\epsilon(t)$ is bounded for every $t>t_0$ and therefore function
$V_\theta$ is well defined. Let us consider derivative
$\dot{V}_\theta$:
\begin{eqnarray}\label{theta_subsystem}
\dot{V}_\theta&=&(\thetavec-\hat{\thetavec})^{T}H(S(\thetavec)-S(\hat{\thetavec}))+
(S(\thetavec)-S(\hat{\thetavec}))^T \times \nonumber
\\
& &H(\thetavec-\hat{\thetavec})-
2(\kappa^2(\xivec)+1)((\thetavec-\hat{\thetavec})^{T}\alphavec^{T}(\xivec)\times
\nonumber
\\
& & \alphavec(\xivec) (\thetavec-\hat{\thetavec}) +
(\thetavec-\hat{\thetavec})^{T}\alphavec^{T}\varevec(t) +
\frac{\|\varevec(t)\|^2}{4})\nonumber
\\
& & =
(\thetavec-\hat{\thetavec})^{T}(S(\thetavec)-S(\hat{\thetavec}))+(S(\thetavec)-S(\hat{\thetavec}))^T\times\nonumber
\\
& & (\thetavec-\hat{\thetavec})- 2 (\kappa^2({\xivec})+1)\times
\nonumber
\\
& &
\|(\thetavec-\hat{\thetavec})^T\alphavec^{T}(\xivec)+0.5\varevec(t)\|^2
\end{eqnarray}
Function $S(\cdot)$ is continuous, therefore, applying Hadamard
lemma we can write the difference
$S(\thetavec)-S(\hat{\thetavec})$ as follows:
$S(\thetavec)-S(\hat{\thetavec})=\int_0^1\frac{\pd
S(\bfz(\lambda))}{\pd \bfz(\lambda)}d\lambda
(\thetavec-\hat{\thetavec})$, $\bfz(\lambda)=\thetavec\lambda +
\hat{\thetavec}(1-\lambda)$. Hence applying Mean Value Theorem we
derive the following $S(\thetavec)-S(\hat{\thetavec})=\frac{\pd
S(\bfz(\lambda'))}{\pd \bfz(\lambda')}(\thetavec-\hat{\thetavec})$
for some $\lambda'\in [0,1]$. The last equation leads to the
following estimation of $\dot{V}_\theta$:
\begin{eqnarray}\label{theta_subsystem_Vdot}
\dot{V}_\theta&
 = &(\thetavec-\hat{\thetavec})^{T}(\frac{\pd S(\bfz(\lambda'))}{\pd
\bfz(\lambda')}^T H+H \frac{\pd S(\bfz(\lambda'))}{\pd
\bfz(\lambda')})(\thetavec-\nonumber \\
& & \hat{\thetavec})\nonumber - 2
(\kappa^2({\xivec})+1)\|(\thetavec-\hat{\thetavec})^T\alphavec^{T}(\xivec)+0.5\varevec(t)\|^2\nonumber
\\
&\leq&-2
(\kappa^2({\xivec})+1)\|(\thetavec-\hat{\thetavec})^T\alphavec^{T}(\xivec)+\nonumber
\\ && 0.5\varevec(t)\|^2\leq 0
\end{eqnarray}
Inequality (\ref{theta_subsystem_Vdot}) ensures that
$(\thetavec-\hat{\thetavec})\in L_\infty$. Taking into account
that for every $\thetavec_0\in\Theta$ solutions
$\thetavec(t,\thetavec_0,t_0)\subset\Omega(\Theta)\subseteq\Theta$
where $\Theta$ is the bounded set, we can conclude that
trajectories $\hat{\thetavec}(t)$ are bounded, i.e.
$\hat{\thetavec}(t)\in L_\infty$. Thus statement 1) is proven.

Let us prove statement 2) of the lemma. Notice that function
$V(\thetavec,\hat{\thetavec},t)$ is non-increasing and bounded
from below. Therefore
$\kappa({\xivec})((\thetavec-\hat{\thetavec})^T\alphavec^{T}(\xivec)+0.5\varevec(t))\in
L_2$. Hence function
$\kappa({\xivec})(\thetavec-\hat{\thetavec})^T\alphavec^{T}(\xivec)$
belongs to $L_2$ as a sum of two functions from $L_2$. The fact
that $\kappa^2(\xivec)+1$ is separated from zero implies that
$(\thetavec-\hat{\thetavec})^T\alphavec^{T}(\xivec)\in L_2$. This
proves statement 2).

Let us show that $\bfx(t)\in L_\infty$ under conditions formulated
in statement 3) of the lemma. Consider derivative
\begin{eqnarray}\label{psi_subsystem_Vdot}
\dpsi&=&\frac{\pd \psi}{\pd \bfx}\bff_0(\bfx) + \frac{\pd
\psi(\bfx)}{\pd
\bfx}\alphavec(\xivec)(\thetavec-\hat{\thetavec})+\frac{\pd
\psi}{\pd \bfx}\varevec(t)\nonumber \\
&=& \frac{\pd \psi}{\pd \bfx}\bff_0(\bfx) + (\frac{\pd
\psi(\bfx)}{\pd \bfx}-\frac{\pd \psi(\xivec)}{\pd
\xivec})(\alphavec(\xivec)(\thetavec-\hat{\thetavec})+\nonumber
\\ & & \varevec(t))+ \frac{\pd \psi(\xivec)}{\pd
\xivec}(\alphavec(\xivec)(\thetavec-\hat{\thetavec})+\varevec(t))
\end{eqnarray}
Notice that  $\psi\in C^1$, $\bfx-\xivec\in L_\infty$ imply that
the norm: $\|\frac{\pd \psi(\bfx)}{\pd \bfx}-\frac{\pd
\psi(\xivec)}{\pd \xivec}\|$ is bounded. Moreover, $\|\frac{\pd
\psi(\xivec)}{\pd \xivec}\|\leq \kappa(\xivec)$. Hence we can
rewrite (\ref{psi_subsystem_Vdot}) as follows:
\begin{equation}\label{psi_subsystem_Vdot1}
\dot{\psi}=\frac{\pd \psi(\bfx)}{\pd \bfx}\bff_0(\bfx)+\mu(t), \
\mu(t)\in L_2
\end{equation}
Function $\beta(\bfx)$ is separated from zero, i.e. $\exists
\delta>0: \ \beta(\bfx)>2\delta\ \forall \bfx\in \Real^n$. Let us
consider the following positive-definite function:
\begin{equation}\label{V_psi}
V_\psi= \int_{0}^{\psi} \varphi(\sigma)d\sigma +
\frac{1}{4\delta}\int_{t}^{\infty} \mu^2(\tau) d\tau
\end{equation}
Taking into account Assumption \ref{asseme:manifold_stability} and
equality (\ref{psi_subsystem_Vdot1}) derivative $\dot{V}_\psi$ can
be estimated as follows:
\begin{equation}
\begin{split}
\dot{V}_\psi&\leq
-\beta(\bfx)\varphi^2(\psi)+\varphi(\psi)\mu(t)-\frac{1}{4\delta}\mu^{2}(t)\nonumber
\\
&\leq -2\delta
\varphi^2(\psi)+\varphi(\psi)\mu(t)-\frac{1}{4\delta}\mu^{2}(t)\nonumber
\\
&= -\delta\varphi^2(\psi)-\delta
(\varphi(\psi)-\frac{1}{2}\mu(t))^2\leq 0\nonumber
\end{split}
\end{equation}
Boundedness of $\bfx$ then follows explicitly from Assumption
 \ref{assume:constraint}. This proves statement 3).

Let us prove that estimate $\hat{\thetavec}(t)$ converges to
$\hat{\thetavec}$ exponentially fast under assumption of
persistent excitation and assuming that $\varevec\equiv 0$.
Consider the following subsystem
\begin{equation}\label{theta_dynamics}
\begin{split}
\dot{\tilde{\thetavec}}&=S(\thetavec)-S(\hat{\thetavec})-H^{-1}(\kappa^2(\xivec)+1)\times
\\
& \alphavec(\xivec)^{T}\alphavec(\xivec)\tilde{\thetavec} =
(\int_0^1\frac{\pd S(\bfz(\lambda))}{\pd \bfz(\lambda)}d\lambda-
\\
&
H^{-1}(\kappa^2(\xivec)+1)\alphavec(\xivec)^{T}\alphavec(\xivec))\tilde{\thetavec}
\end{split}
\end{equation}
where $\tilde{\thetavec}=\thetavec-\hat{\thetavec}$. According to
equations (\ref{system_affine_aux}) system (\ref{theta_dynamics})
describe dynamics of $\hat{\thetavec}(t)-\hat{\thetavec}(t)$.
Solution of (\ref{theta_dynamics}) can be derived in the following
form ${\tilde{\thetavec}}(t)=e^{\int_0^t \frac{\pd
S(\thetavec'(\tau))}{\pd \thetavec'}}d\tau$ $e^{-H^{-1}\int_0^t
(\kappa^2(\xivec(\tau))+1)\alphavec^T(\xivec(\tau))\alphavec(\xivec(\tau))d\tau}$
$\times\tilde{\thetavec}(t_0)$, where
$\thetavec'(\tau)=\thetavec(\tau)\lambda'-\hat{\thetavec}(\tau)(1-\lambda')$
for some $\lambda\in[0,1]$. It follows from Assumption
\ref{asume:internal_model} that the induced matrix norm of
$e^{\int_0^t \frac{\pd S(\thetavec'(\tau))}{\pd \thetavec'}}d\tau$
is bounded, i. e. there exists some positive $D_0>0$ such that
$\|e^{\int_0^t \frac{\pd S(\thetavec'(\tau))}{\pd
\thetavec'}}d\tau\|\leq D_0$ for all $t\geq 0$. On the other hand,
for every $t>T$ there exists integer $ n >0$ such that $t=n T +r,
\ r\in\Real_+ < T$, and the following estimation holds: $
\|e^{-H^{-1}\int_0^t (\kappa^2(\xivec(\tau))+1)
\alphavec^T(\xivec(\tau))\alphavec(\xivec(\tau))d\tau}\|$ $\leq
D_0 \|e^{-H^{-1}\delta I_d n }\|$ $\leq
\|e^{-H^{-1}\frac{\delta}{T} I_d t + I}\|$. Hence we can bound
the norm $\|\tilde{\thetavec}(t)\|$ as follows:
\[
\|\tilde{\thetavec}(t)\|\leq  D_0 \|e^{-H^{-1}\frac{\delta}{T} I_d
t + I}\|\|\tilde{\thetavec}(t_0)\|
\]
{\it The lemma is proven.}

\subsection{Embedding (design of the extension)}\label{proof:part_2}

In this section we show that for the class of systems given by
(\ref{system_affine}) with locally Lipshitz  $\phi_i(\bfx)$:
\begin{equation}\label{phi_def}
\begin{split}
 &  \phi(\bfx): \Real^n\rightarrow \Real^{d\times m}, \nonumber\\
&  \phi(\bfx)=\left( \begin{array}{lll}
                   \phi_{1,1}(\bfx), &  \dots, &  \phi_{1,d}(\bfx)\\
                   \dots, &  \dots, &  \dots \\
                   \phi_{m,1}(\bfx), &  \dots, &  \phi_{m,d}(\bfx)
                   \end{array}
                   \right), \nonumber\\
&  \phi_i(\bfx) =
(\phi_{i,1}(\bfx),\dots,\phi_{i,d}(\bfx))\nonumber
\end{split}
\end{equation}
one can design $C^{1}$-smooth function $\xivec(t)$ such that
$(\alphavec(\bfx)-\alphavec(\xivec))\thetavec(t)$,
$\kappa(\xivec)(\alphavec(\bfx)-\alphavec(\xivec))\thetavec(t)\in
L_2$.

\begin{lem}\label{lemma:extension} Let system (\ref{system_affine})
be given and functions $\phi_i(\bfx)$ defined as in
(\ref{phi_def}) be locally Lipshitz:
\[
\|\phi_i(\bfx)-\phi_i(\xivec)\|\leq
\lambda_i(\bfx,\xivec)\|\bfx-\xivec\|,
\]
where $\lambda(\bfx,\xivec):\Real^n\times \Real^{n}\rightarrow
\Real_{+}$, $\lambda(\bfx,\xivec)$ is locally bounded w.r.t.
$\bfx$, $\xivec$. Let, furthermore, Assumption
\ref{asume:internal_model} hold. Then there exists system
\begin{eqnarray}\label{ext_system}
\dot{\xivec}&=&\bff(\bfx)+ G_u \bfu +  \lambda(\bfx,\xivec)(\bfx-\xivec)+G_u\phi(\bfx)\nuvec\nonumber\\
\dot{\nuvec}&=&S(\nuvec)+H^{-1}(G_u\phi(\bfx))^{T}
(\bfx-\xivec)^{T},
\nonumber \\
\lambda(\bfx,\xivec)&=&1+\sum_{i=1}^{m}\lambda_i^2(\bfx,\xivec)(1+\kappa^2(\xivec))
\end{eqnarray}
such that the following hold:

1)  $\|(\alphavec(\bfx)-\alphavec(\xivec))\thetavec\|\in L_2$,
$\|\kappa(\xivec)(\alphavec(\bfx)-\alphavec(\xivec))\thetavec\|\in
L_2$ for every bounded $\thetavec$;

2) $\bfx\in L_\infty\Rightarrow \xivec\in L_\infty$,
$\lim_{t\rightarrow\infty}\bfx(t)-\xivec(t)=0$
\end{lem}
{\it Proof of Lemma \ref{lemma:extension}.} To prove the lemma it
is enough to consider the following positive definite function
$V_\xi$:
\[
V_\xi=0.5\|\bfx-\xivec\|^2+0.5\|\thetavec-\nuvec\|_{H}^2.
\]
Its time-derivative can be written as follows:
\begin{equation}
\begin{split}
 & \dot{V}_\xi\leq
-\lambda(\bfx,\xivec)\|\bfx-\xivec\|^2+(\bfx-\xivec)^T
G_u\phi(\bfx)(\thetavec-\nuvec)\nonumber
\\
& + (\thetavec-\nuvec)^T(G_u\phi(\bfx))^{T}(\bfx-\xivec)\leq
-\lambda(\bfx,\xivec)\|\bfx-\xivec\|^2\nonumber
\end{split}
\end{equation}
The last inequality implies that
$\lambda_i(\bfx,\xivec)\|\bfx-\xivec\|$,
$\kappa(\xivec)\lambda_i(\bfx,\xivec)\|\bfx-\xivec\| \in L_2$.
Hence
\begin{equation}
\begin{split}
&  \|\phi_i(\bfx)-\phi_i(\xivec)\|\leq
\lambda_i(\bfx,\xivec)\|\bfx-\xivec\|\Rightarrow \nonumber \\
&  \|\phi_i(\bfx)-\phi_i(\xivec)\|,
\kappa(\xivec)\|\phi_i(\bfx)-\phi_i(\xivec)\|\in L_2\nonumber
\end{split}
\end{equation}
Therefore, boundedness of $\thetavec(t)$ and finiteness of the
induced norm $G_u$ ensure that
$\|G_u(\phi_i(\bfx)-\phi_i(\xivec))\thetavec(t)\|$,
$\|\kappa(\xivec)G_u(\phi_i(\bfx)-\phi_i(\xivec))\thetavec(t)\|\in
L_2$.

In order to complete the proof we notice that function $V_\xi$ is
nonincreasing and radially unbounded. This guarantees that
$\xivec$ is bounded as long as $\bfx$ remains bounded. The fact
that $\lambda(\bfx,\xivec)>1$ implies that $\bfx-\xivec\in L_2$.
Under assumptions of the lemma, the right-hand side of the system
is locally bounded. This leads to uniform continuity of
$\|\bfx-\xivec\|^2$, which guarantees that $\lim_{t\rightarrow
\infty}(\bfx-\xivec)=0$. {\it The lemma is proven}

\subsection{Embedding (proof of Theorems \ref{theorem:adaptive_invariant}, \ref{theorem:adaptive_invariant_integral})}\label{proof:part_3}

In this section we provide technical proof of the main results of
our paper.

{\it Proof of Theorem \ref{theorem:adaptive_invariant}.} According
to Lemma \ref{lemma:extension} there exist system
(\ref{ext_system}):
\begin{equation}\label{ext_system1}
\begin{split}
\dot{\xivec}&=\bff(\bfx)+ G_u \bfu +
\lambda(\bfx,\xivec)(\bfx-\xivec)+\\
& G_u\phi(\bfx)\nuvec\\
\dot{\nuvec}&=S(\nuvec)+H^{-1}(G_u\phi(\bfx))^{T}
(\bfx-\xivec)^{T},
\\
\lambda(\bfx,\xivec)&=1+\sum_{i=1}^{m}\lambda_i^2(\bfx,\xivec)(1+\kappa^2(\xivec))
\end{split}
\end{equation}
such that  $\|G_u(\phi(\bfx)-\phi(\xivec))\thetavec\|$,
$\|\kappa(\xivec)G_u(\phi(\bfx)-\phi(\xivec))\thetavec\|\in L_2$
for every bounded $\thetavec(t)$ and trajectory $\bfx(t)$
generated by
\begin{equation}\label{system_controlled}
\dot{\bfx}=\bff(\bfx)+G_u (\phi(\bfx)\thetavec+\bfu); \
\dot{\thetavec}=S(\thetavec)
\end{equation}
Using the notation introduced in the previous subsections:
$\alphavec(\xivec)=G_u\phi(\xivec)$, taking into account that
$\bfu(\bfx,\hat{\thetavec})=\bfu_0(\bfx)-\phi(\xivec)\hat{\thetavec}(t)$,
and denoting
$\varevec(t)=(\alphavec(\bfx)-\alphavec(\xivec))\thetavec(t)$ we
rewrite (\ref{system_controlled}) as follows:
\begin{equation}\label{system_controlled1}
\begin{split}
\dot{\bfx}&=\bff_0(\bfx)+\alphavec(\xivec)(\thetavec-\hat{\thetavec}(t))+\varevec(t)\\
\dot{\thetavec}&=S(\thetavec), \ \varevec(t)\in L_2, \
\kappa(\xivec)\varevec(\xivec)\in L_2
\end{split}
\end{equation}
Taking into account equation (\ref{system_controlled1}) and
expression (\ref{fifo_lin_invar_set}) specifying the function
$\hat{\thetavec}(t)$ we can derive the time-derivative
$\dot{\hat{\thetavec}}$:
\begin{eqnarray}\label{theta_subsystem_eq}
\dot{\hat{\thetavec}}&=&S(\hat{\thetavec})+H^{-1}(\kappa^2(\xivec)+1)\alphavec^{T}(\xivec)(\alphavec(\xivec)(\thetavec-\hat{\thetavec})\nonumber
\\ & & +\varevec(t))
\end{eqnarray}
Then applying Lemma \ref{lemma:virtual} we can conclude that both
$\bfx(t)$ and $\hat{\thetavec}$ are bounded, i.e.
$\bfx(t),\hat{\thetavec}(t)\in L_\infty$. On the other hand,
according to Lemma \ref{lemma:extension}, boundedness of $\bfx$
implies boundedness of $\xivec(t)$. Hence statement 1) of the
theorem is proven.

Notice also that according to Lemma \ref{lemma:extension} the
following holds: $\bfx(t)-\xivec(t)\rightarrow 0$ as
$t\rightarrow\infty$. This fact together with uniform asymptotic
stability of unperturbed system (\ref{theta_subsystem_eq}) (i. e.
when  $\varevec(t)\equiv 0$) imply that
$\hat{\thetavec}(t,\hat{\thetavec}_0,t_0)\rightarrow
\thetavec(t,\thetavec_0,t_0)$ as $t\rightarrow \infty$.This proves
statement 3) of the theorem.

Let us prove that $\bfx(t)\rightarrow \Omega^{\ast}$ as
$t\rightarrow \infty$. In order to do this let us rewrite the
closed-loop system in the following form:
\begin{equation}\label{system_combined}
\begin{split}
\dot{\bfx}&=\bff_0(\bfx)+\alphavec(\xivec)(\thetavec-\hat{\thetavec}(t))+\varevec(t)\\
\dot{\thetavec}&=S(\thetavec) \\
\dot{\hat{\thetavec}}&=
S(\hat{\thetavec})+H^{-1}(\kappa^2(\xivec)+1)\alphavec(\xivec)^{T}\times \\
&  (\alphavec(\xivec)(\thetavec-\hat{\thetavec})+\varevec(t))\\
\dot{\xivec}&=\bff(\bfx)+ G_u \bfu +  \lambda(\bfx,\xivec)(\bfx-\xivec)+ \\
&  G_u\phi(\bfx)\nuvec\\
\dot{\nuvec}&=S(\nuvec)+H^{-1}(G_u\phi(\bfx))^{T}
(\bfx-\xivec)^{T},
 \\
\varevec(t)&= (\alphavec(\bfx)-\alphavec(\xivec))\thetavec, \
\varevec\in L_2 \\
\dot{\epsilon}_0&= - \|\frac{\pd \psi(\bfx)}{\pd
\bfx}(\alphavec(\xivec)(\thetavec-\hat{\thetavec})+\varevec)\|^2
\\
\dot{\epsilon}_1&=
-\|(\alphavec(\bfx)-\alphavec(\xivec))\thetavec\|^{2} \\
\dot{\epsilon}_2&=
-\|\alphavec(\xivec)(\thetavec-\hat{\thetavec})\|^{2}
\end{split}
\end{equation}
I has been shown earlier that trajectories of system
(\ref{system_combined}) are bounded except for the function
$\epsilon_0$. Boundedness of $\epsilon_0(t)$, however, follows
immediately from the fact that $\frac{\pd \psi(\bfx)}{\pd \bfx}$
is bounded and that $\varevec,
\alphavec(\xivec)(\thetavec-\hat{\thetavec}) \in L_2$. Let us
consider the following function:
$V=\int_{0}^\psi\varphi(\sigma)d\sigma +
\|\thetavec-\hat{\thetavec}\|^{2}_{H}+
\frac{1}{4\delta}\epsilon_0(t) + \frac{1}{4}\epsilon_1(t) +
\epsilon_2(t)$. Its time-derivative satisfies the following
inequality: $\dot{V}\leq - \delta \varphi^2(\psi)-
\|\alphavec(\xivec)(\thetavec-\hat{\thetavec})+0.5 \varevec(t)
\|^2$ $-\|\alphavec(\xivec)(\thetavec-\hat{\thetavec})+\varevec(t)
\|^2$ $\leq -\delta \varphi^2(\psi)
-\|\alphavec(\xivec)(\thetavec-\hat{\thetavec})+\varevec(t) \|^2$
Therefore, applying LaSalle invariance principle \cite{LaSalle} we
can conclude that $(\bfx(t),\hat{\thetavec}(t))$ converge (as
$t\rightarrow \infty$) to the largest invariant set in
$\Omega_\psi\times \Omega_\theta$, where
$\Omega_{\psi}=\{\bfx\in\Real^{n}| \ \bfx: \varphi(\psi(\bfx))=0
\}$, and $\Omega_\theta:\{\hat{\thetavec}\in \Real^{d} \ | \
\hat{\thetavec}:
\alphavec(\xivec)(\thetavec-\hat{\thetavec})+\varevec(t)=0\}$. For
the trajectory $\bfx(t)$ this set is defined as the largest
invariant set of system
\begin{equation}\label{system_invariant}
\dot{\bfx}=\bff_0(\bfx)
\end{equation}
under restriction that $\bfx(t)\in\Omega_\psi$.  According to
Assumption \ref{assume:target_invariance} the largest invariant
set of (\ref{system_invariant}) in $\Omega_\psi$ is
$\Omega^{\ast}$. Q.E.D.

{\it Proof of Theorem \ref{theorem:adaptive_invariant_integral}.}
Consider system (\ref{ext_system1}). It follows from Lemma
\ref{lemma:extension} and Assumption \ref{asume:internal_model}
that $G_u(\phi(\bfx)-\phi(\xivec))\thetavec\in L_2$. Then
boundedness of $\hat{\thetavec}(t)$ follows explicitly from the
proof of Theorem \ref{theorem:adaptive_invariant} (let
$\kappa(\xivec)\equiv 0$ in (\ref{theta_subsystem_Vdot})).
Furthermore, Lemma \ref{lemma:virtual} ensures that
$G_u\phi(\xivec)(\thetavec-\hat{\thetavec})\in L_2$. Hence
denoting
$\varevec_0(t)=G_u\phi(\xivec)(\thetavec-\hat{\thetavec})+G_u(\phi(\bfx)-\phi(\xivec))\thetavec$
we obtain that trajectories $\bfx(t)$ in system
(\ref{system_affine}) satisfy the following equation:
\begin{equation}\label{system_controlled2}
\dot{\bfx}=\bff_0(\bfx)+\varevec(t),
\end{equation}
where  $\varevec(t)\in L_2$. System (\ref{system_controlled2}),
however, has finite $L_2\rightarrow L_\infty$ gain and therefore
$\bfx(t)$ is bounded. Therefore, statement 1) of the theorem is
proven. Statement 3) follows explicitly from Lemma
\ref{lemma:virtual}. Let us show that $\bfx(t)\rightarrow
\Omega^{\ast}$ as $t\rightarrow\infty$. In order to do so let us
consider system (\ref{system_combined}) excluding the equation for
$\epsilon_0$. We have already shown that solutions of system
(\ref{system_combined}) are bounded. Define $V=
\|\thetavec-\hat{\thetavec}\|^{2}_{H}+ \frac{1}{4}\epsilon_1(t) +
\epsilon_2(t)$. Its time-derivative satisfies the following
inequality: $\dot{V}\leq
-\|\alphavec(\xivec)(\thetavec-\hat{\thetavec})+\varevec(t) \|^2$
and therefore, applying LaSalle invariance principle
\cite{LaSalle} we obtain that $\bfx(t)\rightarrow \Omega^{\ast}$
as $t\rightarrow\infty$. {\it The theorem is proven.}

\section{Conclusion}

In this paper we have proposed a new framework for adaptive
regulation to invariant sets. The main advantage of our approach
is that we do not require knowledge of the strict Lyapunov
functions for design of the adaptation schemes. Our method also
handles non-equilibrium desired regimes of the system. In addition
it does not assume asymptotic Lyapunov stability of the taget
dynamics.

The number  of the additional equations required for
implementation of our method is  $(n+2 d)$ which compares
favorably with $(n d + d + n)$ in \cite{Panteley_2002}. Though the
conditions we require differ from that of  \cite{Panteley_2002},
we believe that our results naturally complement the existing ones
without too much of additional restrictions.

In the present study we considered linear parameterizations of the
uncertainties. On the other hand, the machinery we use in the
proofs allows to extend the results to nonlinear parameterized
systems \cite{tpt2003_tac,t_fin_forms_arch}. This together with
robustness analysis are currently the the topics of our future
studies.

\bibliography{robust_detect_lin_conf}

\end{document}